\documentclass[11pt,leqno,fleqn]{article}
\usepackage{epsf,amsmath,amssymb,hyperref,srcltx}
\usepackage{amsfonts}
\usepackage{psfrag}
\usepackage{nicefrac}
\usepackage[usenames]{color}
\definecolor{Red}{rgb}{1.0,0.1,0.0}
\newcommand{\Sym}{\text{Sym}}
\newcommand{\dps}{\displaystyle}
\newcommand{\End}{\text{\rm End}}
\newcommand{\kies}[2]{\mbox{${{#1}\choose{#2}}$}}
\newcommand{\dy}[2]{%
\refstepcounter{equation}%
\label{#1}%
\begin{list}{}{
\topsep 5mm
\leftmargin 18mm
\rightmargin 0cm
\itemsep 0mm
\listparindent 0mm
\parsep 0mm
\itemsep 0mm
\labelsep 0mm
\labelwidth 18mm
}%
\item[\rm (\theequation)\hfill]
#2
\end{list}%
}
\newcommand{\dyz}[1]{%
\refstepcounter{equation}%
\begin{list}{}{
\topsep 5mm
\leftmargin 18mm
\rightmargin 0cm
\itemsep 0mm
\listparindent 0mm
\parsep 0mm
\itemsep 0mm
\labelsep 0mm
\labelwidth 18mm
}%
\item[\rm (\theequation)\hfill]
#1
\end{list}%
}
\newcommand{\dyyz}[1]{\dyz{\raggedright$\dps#1$}}
\newcommand{\dyy}[2]{\dy{#1}{\raggedright$\dps#2$}}
\newcommand{\de}[2]{\dy{#1}{\raggedright$\displaystyle #2 $}}
\newcommand{\dez}[1]{\dyz{\raggedright$\displaystyle #1 $}}
\newcounter{bewering}
\newcommand{\prop}[2]{\refstepcounter{bewering}\vspace{4mm}\noindent{\bf Proposition \thebewering.}\label{#1}{\it #2}}
\newcommand{\propz}[1]{\refstepcounter{bewering}\vspace{4mm}\noindent{\bf Proposition \thebewering.}{\it #1}}
\newcounter{sectie}
\newcounter{subsectie}
\newcommand{\sect}[2]{\refstepcounter{sectie}\setcounter{subsectie}{0}
\section*{\boldmath \thesectie. #2}%
\label{#1}}
\newcommand{\sectz}[1]{\refstepcounter{sectie}\setcounter{subsectie}{0}
\section*{\boldmath \thesectie. #1}%
}
\newcommand{\subsect}[2]{\refstepcounter{subsectie}
\subsection*{\boldmath \thesectie.\thesubsectie. #2}%
\label{#1}}
\newcommand{\subsectz}[1]{\refstepcounter{subsectie}
\subsection*{\boldmath \thesectie.\thesubsectie. #1}%
}
\newcommand{\pf}{\vspace{3mm}\noindent{\bf Proof.}\ }
\newcommand{\bx}{\hspace*{\fill} \hbox{\hskip 1pt \vrule width 4pt height 8pt depth 1.5pt \hskip 1pt}

\addvspace{4mm}}
\newcommand{\bxx}{\hspace*{\fill} \hbox{\hskip 1pt \vrule width 4pt height 8pt depth 1.5pt \hskip 1pt}}
\newcommand{\rf}[1]{{\rm (\ref{#1})}}
\newcommand{\AAA}{{\cal A}}
\newcommand{\BB}{{\cal B}}

\newcommand{\EE}{{\cal E}}

\newcommand{\NN}{{\cal N}}

\newcommand{\PP}{{\cal P}}

\newcommand{\kfrac}[2]{\mbox{$\frac{#1}{#2}$}}

\newcommand{\oC}{{\mathbb{C}}}

\newcommand{\oR}{{\mathbb{R}}}

\newcommand{\oZ}{{\mathbb{Z}}}
\begin{document}
\begin{center}
{\LARGE\bf Semidefinite code bounds based on quadruple distances

}

\end{center}
\vspace{1mm}
\begin{center}
{\large
\hspace{10mm}
Dion C. Gijswijt\footnote{ CWI and Department of Mathematics,
Leiden University},
Hans D. Mittelmann\footnote{ School of Mathematical and Statistical
Sciences, Arizona State University},
and
Alexander Schrijver\footnote{ CWI and Department of Mathematics, University of Amsterdam.
Mailing address: CWI, Science Park 123, 1098 XG Amsterdam,
The Netherlands.
Email: lex@cwi.nl.}}

\end{center}

\noindent
{\small{\bf Abstract.}
Let $A(n,d)$ be the maximum number of $0,1$ words of length $n$, any two
having Hamming distance at least $d$.
We prove $A(20,8)=256$, which implies that the quadruply shortened Golay code
is optimal.
Moreover, we show
$A(18,6)\leq 673$,
$A(19,6)\leq 1237$,
$A(20,6)\leq 2279$,
$A(23,6)\leq 13674$,
$A(19,8)\leq 135$,
$A(25,8)\leq 5421$,
$A(26,8)\leq 9275$,
$A(21,10)\leq 47$,
$A(22,10)\leq 84$,
$A(24,10)\leq 268$,
$A(25,10)\leq 466$,
$A(26,10)\leq 836$,
$A(27,10)\leq 1585$,
$A(25,12)\leq 55$,
and
$A(26,12)\leq 96$.

The method is based on the positive semidefiniteness of matrices
derived from quadruples of words.
This can be put as constraint in a semidefinite program, whose optimum
value is an upper bound for $A(n,d)$.
The order of the matrices involved is huge.
However, the semidefinite program is highly symmetric, by which its
feasible region can be restricted to the algebra of matrices invariant
under this symmetry.
By block diagonalizing this algebra, the order of the matrices will be
reduced so as to make the program solvable with semidefinite programming
software in the above range of values of $n$ and $d$.

}


\sectz{Introduction}

For any $k$, we will identify elements $\{0,1\}^k$ with
$0,1$ words of length $k$.
The {\em Hamming distance} $d_H(v,w)$ between two words $v,w$
is the number of $i$ with $v_i\neq w_i$.

Throughout we denote
\dez{
N:=\{0,1\}^n.
}
A {\em code} of {\em length} $n$ is any subset $C$ of $N$.
The {\em minimum distance} of a code $C$ is the minimum Hamming distance
between any two distinct elements of $C$.
Then $A(n,d)$ denotes the maximum size (= cardinality) of a code of
length $n$ with minimum distance at least $d$.

Computing $A(n,d)$ and finding upper and lower bounds for it have been
long-time focuses in combinatorial coding theory
(cf.\ MacWilliams and Sloane [11]).
Classical is Delsarte's bound
[3].
Its value can be described as the maximum $A_2(n,d)$ of $\sum_{u,v\in N}X_{u,v}$,
where $X$ is a symmetric, nonnegative, positive semidefinite
$N\times N$ matrix with trace 1 and with $X_{u,v}=0$ if $u,v\in N$ are
distinct and have distance less than $d$.
Then $A(n,d)\leq A_2(n,d)$, since for any nonempty code $C$ of minimum
distance at least $d$, the matrix $X$ with $X_{u,v}=|C|^{-1}$ if
$u,v\in C$ and $X_{u,v}=0$ otherwise, is
a feasible solution with objective value $|C|$.

This is the {\em analytic} definition of the Delsarte bound
(in the vein of Lov\'asz [9], cf.\
[12],
[16]).
It is a semidefinite programming problem, but of huge dimensions ($2^n$),
which makes it hard to compute in this form.

However, the problem is highly symmetric.
The group $G$ of distance preserving permutations of the set $N$ acts on
the set of optimum solutions: if $(X_{u,v})$ is an optimum solution,
then also $(X_{\pi(u),\pi(v)})$ is an optimum solution for any
$\pi\in G$.
Hence, as the set of optimum solutions is convex, by averaging we obtain
a $G$-invariant optimum solution $X$.
That is, $X_{\pi(u),\pi(v)}=X_{u,v}$ for all $u,v$ and all $\pi\in G$.
So $X_{u,v}$ depends only on the Hamming distance of $u$ and $v$,
hence there are in fact at most $n+1$ variables.
Since (in this case) the algebra of $G$-invariant matrices is commutative,
it implies that there is a unitary matrix $U$ such
that $U^*XU$ is a diagonal matrix for each $G$-invariant $X$.
It reduces the semidefinite constraints of order $2^n$ to $2^n$ linear
constraints, namely the nonnegativity of the diagonal elements.
As the space of $G$-invariant matrices is $n+1$-dimensional, there are
in fact only $n+1$ different linear constraints, hence it reduces to a small
linear programming problem.

So the Delsarte bound is initially a huge semidefinite program in
variables associated with pairs and singletons of words in
$N$, that can be reduced to a small linear program, with
a small number of variables.
In [17] this method was generalized to semidefinite
programs in variables associated with sets of words of size
at most 3.
In that case, the programs can be reduced by block diagonalization to a small
semidefinite program, with a small number of variables.
A reduction to a {\em linear} program does not work here, as in this case
the corresponding algebra is not commutative.
This however is not a real bottleneck, as like for linear programming
there are efficient (polynomial-time) algorithms for semidefinite programming.

In the present paper we extend this method to quadruples of words.
Again, by a block diagonalization, the order of the size of the semidefinite
programs is reduced from exponential size to polynomial size.
We will give a more precise description of the method
in Section \ref{11fe10a}.

The reduced semidefinite programs still tend to get rather large,
but yet for $n$ up to $28$ and several values of $d$, we were able
to solve the associated semidefinite programming up to (more than)
enough precision, using the semidefinite programming algorithm SDPA.
It gives the new upper bounds $A_4(n,d)$ for $A(n,d)$ displayed in Table 1.
One exact value follows, namely $A(20,8)=256$.
It means that the quadruply shortened Golay code is optimum.
In the table we give also the values of the new bound where it does not improve the
currently best known bound, as in many of such cases the new bound confirms or
is very close to this best known bound.

Since $A(n,d)=A(n+1,d+1)$ if $d$ is odd, we can restrict ourselves
to $d$ even.
We refer to the websites maintained by
Erik Agrell [1] and Andries Brouwer [2]
for more background on the known upper and lower bounds displayed in the table.

\begin{table}[h!]
\begin{center}
\begin{tabular}{|r|r|r|r|r||r@{.}l|}
\hline
&&known&known&\textcolor{Red}{new}&\multicolumn{2}{|c|}{}\\
$n$&$d$&lower&upper&\textcolor{Red}{upper}&\multicolumn{2}{|c|}{$A_4(n,d)$}\\
&&bound&bound&\textcolor{Red}{bound}&\multicolumn{2}{|c|}{}\\
\hline
\hline
17&4&2720&3276&&
3276&800\\
18&4&5312&6552&&
6553&600\\
19&4&10496&13104&&
13107&200\\
20&4&20480&26168&&
26214&400\\
21&4&36864&43688&&
43690&667\\
\hline
17&6&256&340&&
351&506\\
\textcolor{Red}{18}&\textcolor{Red}{6}&512&680&\textcolor{red}{673}&
673&005\\
\textcolor{Red}{19}&\textcolor{Red}{6}&1024&1280&\textcolor{red}{1237}&
1237&939\\
\textcolor{Red}{20}&\textcolor{Red}{6}&2048&2372&\textcolor{red}{2279}&
2279&758\\
21&6&2560&4096&&
4096&000\\
22&6&4096&6941&&
6943&696\\
\textcolor{Red}{23}&\textcolor{Red}{6}&8192&13766&\textcolor{red}{13674}&
13674&962\\
\hline
17&8&36&37&&
38&192\\
18&8&64&72&&
72&998\\
\textcolor{Red}{19}&\textcolor{Red}{8}&128&142&\textcolor{red}{135}&
135&710\\
\textcolor{Red}{20}&\textcolor{Red}{8}&\textcolor{red}{256}&274&\textcolor{red}{256}&
256&000\\
\textcolor{Red}{25}&\textcolor{Red}{8}&4096&5477&\textcolor{red}{5421}&
5421&499\\
\textcolor{Red}{26}&\textcolor{Red}{8}&4096&9672&\textcolor{red}{9275}&
9275&544\\
\hline
\textcolor{Red}{21}&\textcolor{Red}{10}&42&48&\textcolor{red}{47}&
47&007\\
\textcolor{Red}{22}&\textcolor{Red}{10}&64&87&\textcolor{red}{84}&
84&421\\
23&10&80&150&&
151&324\\
\textcolor{Red}{24}&\textcolor{Red}{10}&128&280&\textcolor{red}{268}&
268&812\\
\textcolor{Red}{25}&\textcolor{Red}{10}&192&503&\textcolor{red}{466}&
466&809\\
\textcolor{Red}{26}&\textcolor{Red}{10}&384&886&\textcolor{red}{836}&
836&669\\
\textcolor{Red}{27}&\textcolor{Red}{10}&512&1764&\textcolor{red}{1585}&
1585&071\\
\hline
\textcolor{Red}{25}&\textcolor{Red}{12}&52&56&\textcolor{red}{55}&
55&595\\
\textcolor{Red}{26}&\textcolor{Red}{12}&64&98&\textcolor{red}{96}&
96&892\\
27&12&128&169&&
170&667\\
28&12&178&288&&
288&001\\
\hline
\end{tabular}

\bigskip
{\bf Table 1.} Bounds for $A(n,d)$
\end{center}
\end{table}

In the computations, the accuracy of the standard double precision version of SDPA
(already considered in the comparison [13]) was insufficient
for several of the cases solved here.
The semidefinite programs generated appear to have rather thin feasible
regions so that SDPA and the other high-quality but double precision codes
terminate prematurely with large infeasibilities.
We have used the multiple precision versions of SDPA developed by
M. Nakata for quantum chemistry computations in [15].
Fortunately, the quadruple precision version was sufficient in
all cases tabulated.
The even higher precision versions as also implemented
by the second author for interactive use at the NEOS Server [14]
would have needed excessive computing times.
Still, the times needed in Table 1 below varied from a few hours for the
small cases to 1$\nicefrac12$ days for $A_4(20,4)$,
13 days for $A_4(23,6)$, 22 days for $A_4(25,8)$,
30 days for $A_4(27,10)$,
and 43 days for $A_4(26,8)$.

The approach outlined above of course suggests a hierarchy of upper bounds
by considering sets of words of size at most $k$, for $k=2,3,4,\ldots$.
This connects to hierarchies of bounds for $0,1$ programming problems
developed by
Lasserre [7],
Laurent [8],
Lov\'asz and Schrijver [10],
and Sherali and Adams [19].
The novelty of the present paper lies in exploiting the symmetry and
giving an explicit block
diagonalization that will enable us to calculate the bounds.

In fact, the relevance of the present paper might be three-fold.
First, it may lie in coding and design theory, as we give new upper bounds
for codes and show that the quadruply shortened Golay code is optimal.
Second, the results may be of interest for algebraic combinatorics
(representations of the symmetric group and extensions),
as we give an explicit block diagonalization of the centralizer
algebra of groups acting on pairs of words from $N$.
Third, the relevance may come from semidefinite programming theory and
practice, by exploiting symmetry and reducing sizes of programs,
and by gaining insight into the border of what is possible with
current-state semidefinite programming software, both as to problem size,
precision, and computing time.

We do not give explicitly all formulas in our description of the method,
as they are sometimes quite involved, rather it may serve as a manual
to obtain an explicit implementation, which should be straightforward
to derive.

\sect{11fe10a}{The bound $A_k(n,d)$}

For any $n,d,k\in\oZ_+$, we define the number
$A_k(n,d)$ as follows.
Let $\NN$ be the collection of codes $S\subseteq N$ of minimum
distance at least $d$.
(Recall that $N=\{0,1\}^n$.)
For any $t$, let $\NN_t$ be the collection of $S\in\NN$ with $|S|\leq t$.

For $S\in\NN_k$, define
\dez{
\NN(S):=\{S'\in\NN\mid S\subseteq S', |S|+2|S'\setminus S|\leq k\}.
}
The rationale of this definition is that
$|S'\cup S''|\leq k$ for all $S',S''\in\NN(S)$.

For $x:\NN_k\to\oR$ and $S\in\NN_k$,
let $M_S(x)$ be the $\NN(S)\times\NN(S)$ matrix
given by
\dez{
M_S(x)_{S',S''}:=
\begin{cases}
x(S'\cup S'')&\text{ if $S'\cup S''\in\NN$,}\\
0&\text{ otherwise,}
\end{cases}
}
for $S',S''\in\NN(S)$.
Define
\dyy{2fe10a}{
A_k(n,d):=
\max\{\sum_{v\in N}x(\{v\})\mid x:\NN_k\to\oR,
x(\emptyset)=1,
M_S(x)$ positive semidefinite for each $S\in\NN_k\}.
}
Note that, as $x(S)$ occurs on the diagonal of $M_S(x)$,
$x$ has nonnegative values only.

\prop{4fe10a}{
$A(n,d)\leq A_k(n,d)$.
}

\pf
Let $C$ be a maximum-size code of length $n$ and minimum distance
at least $d$.
Define $x(S):=1$ if $S\subseteq C$ and $x(S):=0$ otherwise.
Then $M_S(x)$ is positive semidefinite for each $S\in\NN_k$.
Moreover, $A(n,d)=|C|=\sum_{v\in N}x(\{v\})$.
\bx

The upper bound $A_2(n,d)$ can be proved to be equal to the Delsarte bound
[3] (see [5]).
The bound given in [17] is a slight sharpening
of $A_3(n,d)$.

Now to make the problem computationally tractable, let again $G$ be the
group of permutations of $N$ that maintain distances.
Then, if $x$ is an optimum solution of \rf{2fe10a} and $\pi\in G$,
also $x^{\pi}$ is an optimum solution.
(We refer to Section \ref{9fe10b}.\ref{13fe10g} for notation.)
As the feasible region in \rf{2fe10a} is convex, by averaging over
all $\pi\in G$ we obtain a $G$-invariant optimum solution.
So we can reduce the feasible region to those
$x$ that are $G$-invariant.
Then $M_S(x)$ is $G_S$-invariant, where $G_S$ is the $G$-stabilizer
of $S$ (= set of $\pi\in G$ with $\pi(S)=S$).
This allows us to block diagonalize $M_S(x)$,
and to make the problems tractable for larger $n$.
In the coming sections we will discuss how to obtain an explicit block
diagonalization.

It is of interest to remark that the equality $A(20,8)=256$ in fact
follows if we take $k=4$ and require in \rf{2fe10a} only that $M_S(x)$ is
positive semidefinite for all $S$ with $|S|=0$ or $|S|=4$.

%

An observation useful to note (but not used in the sequel) is the
following.
A well-known relation is $A(n+1,d)\leq 2A(n,d)$.
The same relation holds for $A_k(n,d)$:

\propz{
For all $n,d$: $A_k(n+1,d)\leq 2A_k(n,d)$.
}

\pf
Let $x$ attain the maximum \rf{2fe10a} for $A_k(n+1,d)$.
For each $S\subseteq\{0,1\}^n$, let $S':=\{w0\mid w\in S\}$ and
$S'':=\{w1\mid w\in S\}$.
Define $x'(S):=x(S')$ and $x''(S):=x(S'')$ for all $S\in\NN$.
Then $x'$ and $x''$ are feasible solutions of \rf{2fe10a} for
$A_k(n,d)$.
Moreover $\sum_{v\in\{0,1\}^n}(x'(\{v\})+x''(\{v\}))=\sum_{w\in\{0,1\}^{n+1}}x(\{v\})$.
Thus $2A_k(n,d)\geq A_k(n+1,d)$.
\bx

This implies, using $A(20,8)=256$ and $A(24,8)=4096$ (the extended
Golay code), that
$A_4(21,8)=512$,
$A_4(22,8)=1024$,
$A_4(23,8)=2048$, and
$A_4(24,8)=4096$.
We did not display these values in the table, and we do not need to solve the
corresponding semidefinite programming problems.

From now on we will fix $k=4$, and we will use the name $k$ for other
purposes.
In Section \ref{12me10b} we will discuss how to find a further reduction
by considering words of even weights only, which is enough to
obtain the bounds.

\sect{9fe10b}{Preliminaries}

In this section we fix some notation and recall a few basic facts.
Underlying mathematical disciplines are representation theory and
C$*$-algebra, but because the potential readership of this paper
might possess diverse background, we give a brief elementary exposition.

\subsectz{Notation}

We denote
\dyz{
$[s,t]:=\{i\in\oZ\mid s\leq i\leq t\}$ and
$[s,t]_{\text{even}}:=\{i\in\oZ\mid s\leq i\leq t; i$ even$\}$.
}

Throughout this paper, $P$, $T$, and $N$ denote the sets of ordered pairs,
triples, and $n$-tuples of elements of $\{0,1\}$, i.e.,
\dyz{
$P:=\{0,1\}^2$, $T:=\{0,1\}^3$, and $N:=\{0,1\}^n$.
}
As mentioned, we identify elements of $\{0,1\}^t$ with $0,1$
words of length $t$.
We will view $\{0,1\}$ as the field of two elements and add
elements of $P$, $T$, and $N$ modulo 2.

For $\alpha:P\to\oZ_+$, we denote
\dez{
i_{\alpha}:=\alpha(10)+\alpha(11),
j_{\alpha}:=\alpha(01)+\alpha(11),
n_{\alpha}:=\alpha(00)+\alpha(10)+\alpha(01)+\alpha(11).
}

For any finite set $V$ and $n\in\oZ_+$, let
\de{13fe10b}{
\Lambda_V^n:=\{\lambda:V\to\oZ_+\mid\sum_{v\in V}\lambda(v)=n\}.
}
For any $\lambda\in\Lambda_V^n$, let
\dyy{13fe10c}{
\Omega_{\lambda}:=\{\rho:\{1,\ldots,n\}\to V\mid
|\rho^{-1}(v)|=\lambda(v)$ for each $v\in V\}.
}
So $\{\Omega_{\lambda}\mid\lambda\in\Lambda_V^n\}$ is the collection of
orbits on $V^n$ under the natural action of the symmetric group
$S_n$ on $V^n$ (cf.\ Section \ref{9fe10b}.\ref{13fe10g}).

Throughout, $G$ denotes the group of distance preserving permutations
of $N$.
The group consists of all permutations of coordinates followed by
swapping 0 and 1 in a subset of the coordinates.

\subsect{13fe10g}{Group actions}

An {\em action} of a group $H$ {\em on} a set $Z$ is a group homomorphism
from $H$ into the group of permutations of $Z$.
One then says that $H$ {\em acts on} $Z$.
An action of $H$ on $Z$ induces in a natural way actions of $H$ on
derived sets like $Z\times Z$, $\PP(Z)$, $\{0,1\}^Z$, and $\oC^Z$.

If $\pi\in H$ and $z\in Z$, then $z^{\pi}$ denotes the image of $z$
under the permutation associated with $\pi$.
If $H$ acts on $Z$, an element $z\in Z$ is called {\em $H$-invariant}
if $z^{\pi}=z$ for each $\pi\in H$.
The set of $H$-invariant elements of $Z$ is denoted by $Z^H$.

A function $\phi:Z\to Z$ is {\em $H$-equivariant} if
$\phi(z^{\pi})=\phi(z)^{\pi}$ for each $z\in Z$ and each $\pi\in H$.
If $Z$ is a vector space, the collection of $H$-equivariant
endomorphisms $Z\to Z$ is denoted by $\End_H(Z)$.
It is called the {\em centralizer algebra} of the action of $H$ on $Z$.

If $Z$ is a finite set and $H$ acts on $Z$, then there is a natural
isomorphism
\de{13fe10a}{
\End_H(\oC^Z)\cong (\oC^{Z\times Z})^H.
}
If $Z$ is a linear space, the symmetric group $S_n$ acts
naturally on the $n$-th tensor power $Z^{\otimes n}$.
As usual, we denote the subspace of symmetric tensors by
\dez{
\Sym^n(Z):=(Z^{\otimes n})^{S_n}.
}

\subsectz{Matrix $*$-algebras}

A {\em matrix $*$-algebra} is a set of matrices (all of the same order)
that is a $\oC$-linear space and is closed under multiplication and
under taking conjugate transpose ($X\mapsto X^*)$.
If a group $H$ acts on a finite set $Z$, then
$(\oC^{Z\times Z})^H$ is a matrix $*$-algebra.

If $\AAA$ and $\BB$ are matrix $*$-algebras,
a function $\phi:\AAA\to\BB$ is an {\em algebra $*$-homomorphism} if
$\phi$ is linear and maintains multiplication and taking conjugate
transpose.
It is an {\em algebra $*$-isomorphism} if $\phi$ is moreover a bijection.

If $\phi:\AAA\to\BB$ is an algebra $*$-homomorphism and $A\in\AAA$ is
positive semidefinite, then also $\phi(A)$ is positive semidefinite.
(This follows from the fact that any matrix $X$ is positive semidefinite
if and only if $X=YY^*$ for some linear combination $Y$ of $X,X^2,\ldots$.
This last statement can be proved by diagonalizing $X$.)

The sets $\oC^{m\times m}$, for $m\in\oZ_+$, are the
{\em full matrix $*$-algebras}.
An algebra $*$-isomorphism $\AAA\to\BB$ is called a
{\em full block diagonalization} of $\AAA$ if $\BB$
is a direct sum of full matrix $*$-algebras.

Each matrix $*$-algebra has a full block diagonalization --- we need
them explicitly in order to perform the calculations for determining
$A_4(n,d)$.
(A full block diagonalization is in fact unique, up to obvious
transformations: reordering the terms in the sum, and resetting
$X\mapsto U^*XU$, for some fixed unitary matrix $U$, applied to some
full matrix $*$-algebra.)

\subsect{8fe10a}{Actions of $S_2$}

Let $Z$ be a finite set on which the symmetric group $S_2$ acts.
This action induces an action of $S_2$ on $\oC^Z$.
For $\pm\in\{+,-\}$, let
$L_{\pm}:=\{x\in\oR^Z\mid x^{\sigma}=\pm x\}$,
where $\sigma$ is the non-identity element of $S_2$.
Then $L_+$ and $L_-$ are the eigenspaces of $\sigma$.

Let $U_{\pm}$ be a matrix whose columns form an orthonormal basis
of $L_{\pm}$.
The matrices $U_{\pm}$ are easily obtained from the $S_2$-orbits on $Z$.
Then the matrix $[U_+~~U_-]$ is unitary.
Moreover, $U_+^*XU_-=0$ for each $X$ in
$(\oC^{Z\times Z})^{S_2}$.
As $L_+$ and $L_-$ are the eigenspaces of $\sigma$,
the function $X\mapsto U_+^*XU_+\oplus U_-^*XU_-$ defines
a full block diagonalization of $(\oC^{Z\times Z})^{S_2}$.

\subsect{8fe10b}{Fully block diagonalizing $\Sym^n(\oC^{2\times 2})$}

We describe a full block diagonalization
\dez{
\chi:\Sym^n(\oC^{2\times 2})\to
\bigoplus_{k=0}^{\lfloor\frac12n\rfloor}
\oC^{[k,n-k]\times[k,n-k]},
}
as can be derived from the work of Dunkl [4]
(cf.\ Vallentin [20], Schrijver [17]).

To this end, let, for any $\alpha\in\Lambda_P^n$,
\de{13fe10d}{
M_{\alpha}:=\sum_{\rho\in\Omega_{\alpha}}
\bigotimes_{i=1}^nE_{\rho(i)}
\in
\Sym^n(\oC^{2\times 2}).
}
Here, for $c=(c_1,c_2)\in P$,
$E_c$ denotes the $\{0,1\}\times\{0,1\}$ matrix
with $1$ in position $c_1,c_2$ and 0 elsewhere.
Then $\{M_{\alpha}\mid\alpha\in\Lambda_P^n\}$ is a basis of
$\Sym^n(\oC^{2\times 2})$.
(Throughout, we identify $\oC^{2\times 2}$ with
$\oC^{\{0,1\}\times\{0,1\}}$.)

For any $\alpha:P\to\oZ_+$ and $k\in\oZ_+$, define the following
number:
\dyy{10fe10a}{
\gamma_{\alpha,k}:=
\sum_{p=0}^k(-1)^p
\kies{k}{p}
\big(
\kies{\alpha(01)+\alpha(00)-k}{\alpha(01)-p}
\kies{\alpha(01)+\alpha(11)-k}{\alpha(01)-p}
\kies{\alpha(10)+\alpha(00)-k}{\alpha(10)-p}
\kies{\alpha(10)+\alpha(11)-k}{\alpha(10)-p}
\big)^{1/2},
}
and the following $[k,n-k]\times[k,n-k]$ matrix $\Gamma_{\alpha,k}$,
where $n:=n_{\alpha}$:
\dyy{22ja10a}{
(\Gamma_{\alpha,k})_{i,j}:=
\begin{cases}
\gamma_{\alpha,k}&\text{ if $i=i_{\alpha}$ and
$j=j_{\alpha}$,}\\
0&\text{ otherwise,}
\end{cases}
}
for $i,j\in[k,n-k]$.
Now $\chi$ is given by
\dyyz{
\chi(M_{\alpha})=\bigoplus_{k=0}^{\lfloor\frac12n\rfloor}\Gamma_{\alpha,k}
}
for $\alpha\in\Lambda_P^n$.

\sect{9fe10a}{Fully block diagonalizing the matrices $M_S(x)$}

Recall that $\NN_t$ consists of all codes of length $n$,
of minimum distance at least $d$, and size at most $t$.
Let $x:\NN_4\to\oR$ be $G$-invariant.
For any $S\in\NN_4$ and any $\pi\in G$,
we have that $M_S(x)$ is equal to $M_{\pi(S)}(x)$, up to 
renaming row and column indices.
So we need to check positive semidefiniteness of $M_S(x)$ for only one set
$S$ from any $G$-orbit on $\NN_4$.
Moreover, $M_S(x)$ belongs to the centralizer algebra of the $G$-stabilizer
$G_S$ of $S$ (the set of all $\pi\in G$ with $\pi(S)=S$).

Consider any $S\in\NN_4$.
If $|S|=4$, then $M_S(x)$ is a $1\times 1$ matrix with entry $x(S)$.
So $M_S(x)$ is positive semidefinite if and only if $x(S)\geq 0$.

Moreover, if $|S|$ is odd and $|S|\leq 3$, then $M_S(x)$ is a principal
submatrix of $M_R(x)$, where $R$ is any subset of $S$ with $|R|=|S|-1$.
(This because if $S'\supseteq S$ and $|S|+2|S'\setminus S|\leq 4$, then
$|R|+2|S'\setminus R|\leq 4$.)

Concluding it remains to consider checking positive semidefiniteness of
$M_S(x)$ only for $S=\emptyset$ and for one element from each
$G$-orbit of codes $S$ in $\NN$ with $|S|=2$.
We first consider the case $|S|=2$.

Note that the $G$-orbit of any $S\in\NN$ with $|S|=2$
is determined by the distance $m$ between the two elements of $S$.
Hence we can assume $S:=\{{\mathbf 0},u\}$,
where $u$ is the element of $N$ with precisely $m$ 1's,
in positions $1,\ldots,m$.
Then there is a one-to-one relation between
\dez{
N':=\{v\in N\mid d_H({\mathbf 0},v)\in[d,n]\text{ and }d_H(u,v)\in\{0\}\cup[d,n]\}
}
and $\NN(S)$, given by $v\mapsto S\cup\{v\}$.

For any $v,w\in N$, let $\rho_{v,w}:\{1,\ldots,n\}\to P$ be defined by
$\rho_{v,w}(i):=(v_i,w_i)$ for $i=1,\ldots,n$.
This gives an embedding
\dez{
\Phi:
\End_{G_S}(\oC^{\NN(S)})\to
\Sym^m(\oC^{2\times 2})\otimes\Sym^{n-m}(\oC^{2\times 2})
}
defined by
\dez{
\Phi(X):=\sum_{v,w\in N'}
X_{S\cup\{v\},S\cup\{w\}}
M_{\rho_{v,w}}
}
for $X\in\End_{G_S}(\oC^{\NN(S)})$.

The image of $\Phi$ is equal to the linear hull of those
$M_{\alpha}\otimes M_{\beta}$ with $\alpha\in\Lambda_P^m$
and $\beta\in\Lambda_P^{n-m}$ such that
$i_{\alpha}+i_{\beta}\in [d,n]$,
$j_{\alpha}+j_{\beta}\in [d,n]$,
$m-i_{\alpha}+i_{\beta}\in \{0\}\cup[d,n]$,
$m-j_{\alpha}+j_{\beta}\in \{0\}\cup[d,n]$.

With the full block diagonalization \rf{22ja10a} this gives
that the image is equal to the direct sum over $k,l$
of the linear hull of the submatrices of
$\Gamma_{\alpha,k}\otimes \Gamma_{\beta,l}$ induced by the rows
and columns indexed by $(i,i')$ with $i+i'\in [d,n]$ and
$m-i+i'\in \{0\}\cup[d,n]$.

The stabilizer $G_S$ contains a further symmetry, replacing any $c\in N$
by $c+u$ (mod 2).
This leaves $S=\{{\mathbf 0},u\}$ invariant.
It means an action of $S_2$, and the corresponding reduction can be obtained with
the method of Section \ref{9fe10b}.\ref{8fe10a}.

\sect{17me10a}{Fully block diagonalizing $M_{\emptyset}(x)$}

In this section we consider $S=\emptyset$.
Then $\NN(S)=\NN_2$, which is the set of all codes of length $n$,
minimum distance at least $d$, and size at most $2$.
We are out for a full block diagonalization of the centralizer algebra
$\End_G(\oC^{\NN_2})$ of the action of $G$ on $\oC^{\NN_2}$.
This will be obtained in a number of steps.

\subsect{17me10b}{The algebra $\AAA$}

We first consider an algebra $\AAA$ consisting of (essentially) $4\times 4$
matrices.
For any $c\in P=\{0,1\}^2$, let $\overline c:=c+(1,1)$ (mod 2).
Let $\AAA$ be the centralizer algebra of the action of $S_2$ on $P$
generated by $c\mapsto \overline c$ on $c\in P$.
We can find a full block diagonalization with the method of Section
\ref{9fe10b}.\ref{8fe10a}.
We need it explicitly.
Note that
\dyyz{
\AAA=\{A\in\oC^{P\times P}\mid A_{\overline c,\overline d}=A_{c,d}$
for all $c,d\in P\}
}
and that $\AAA$ is a matrix $*$-algebra of dimension 8.

For $c,d\in P$, let $E_{c,d}$ be the $0,1$ matrix in $\oC^{P\times P}$
with precisely one 1, in position $(c,d)$.
Recall $T=\{0,1\}^3$, and define for $t\in T$:
\dez{
B_t:=E_{c,d}+E_{\overline c,\overline d},
}
where $(c,d)$ is any of the two pairs in $P^2$ satisfying
\dyy{15ja10a}{
c_1+c_2=t_1, d_1+d_2=t_2, c_2+d_2=t_3.
}
Then $\{B_t\mid t\in T\}$ is a basis of $\AAA$.

For $i\in\{0,1\}$, let $U_i\in\oC^{P\times\{0,1\}}$ be defined by
\dez{
(U_i)_{c,a}=\kfrac12\sqrt{2}(-1)^{ic_2}\delta_{a,c_1+c_2}
}
for $c\in P$ and $a\in\{0,1\}$.
Then one directly checks that the matrix $U:=[U_0~~U_1]$ is
unitary, i.e., $U^*U=I$.
Moreover, for all $c,d\in P$ and $a,b,i,j\in\{0,1\}$ we have
\dyyz{
(U_i^*E_{c,d}U_j)_{a,b}=
(U_i)_{c,a}(U_j)_{d,b}=
\kfrac12(-1)^{ic_2+jd_2}\delta_{a,c_1+c_2}\delta_{b,d_1+d_2}.
}
Hence, if $t\in T$ and $c,d$ satisfy \rf{15ja10a}, then
\dyy{4no09a}{
(U_i^*B_tU_j)_{a,b}=
\kfrac12((-1)^{ic_2+jd_2}+(-1)^{ic_2+jd_2+i+j})
\delta_{a,c_1+c_2}\delta_{b,d_1+d_2}
=
\kfrac12(-1)^{ic_2+jd_2}(1+(-1)^{i+j})
\delta_{a,c_1+c_2}\delta_{b,d_1+d_2}
=
(-1)^{it_3}
\delta_{i,j}\delta_{a,t_1}\delta_{b,t_2}.
}
So $U_0^*\AAA U_1=0$, and hence, as $\dim\AAA=8$,
$U^*\AAA U$ gives a full block diagonalization of $\AAA$.
Moreover
\dyy{15ja10c}{
U_i^*B_tU_i=(-1)^{it_3}E_{t_1,t_2}.
}

\subsectz{The algebra $\Sym^n(\AAA)$}

Our next step is to find a full block diagonalization
of $\End_G(\oC^{N^2})$, where $N^2$ is (as usual)
the collection of ordered pairs from $N=\{0,1\}^n$.

There is a natural algebra isomorphism
\dez{
\End_{G}(\oC^{N^2})\to\Sym^n(\AAA)
}
by the natural isomorphisms
\de{6ma10a}{
\oC^{(\{0,1\}^n)^2}
\cong
\oC^{(\{0,1\}^2)^n}
\cong
(\oC^{\{0,1\}^2})^{\otimes n},
}
using the fact that $G$ consists of all permutations
of $N$ given by a permutation of the indices in $\{0,\ldots,n\}$
followed by swapping 0 and 1 on a subset of it.

Let $U_0$ and $U_1$ the $P\times\{0,1\}$ matrices given
in Section \ref{17me10a}.\ref{17me10b}.
Define
\dyy{22ja10b}{
\phi:\Sym^n(\AAA)\to
\bigoplus_{m=0}^n
\Sym^m(\oC^{\{0,1\}\times\{0,1\}})
\otimes
\Sym^{n-m}(\oC^{\{0,1\}\times\{0,1\}})
}
by
\dez{
\phi(A):=\bigoplus_{m=0}^n
(U_0^{\otimes m}\otimes U_1^{\otimes n-m})^*
A
(U_0^{\otimes m}\otimes U_1^{\otimes n-m})
}
for $A\in\Sym^n(\AAA)$.
Trivially, $\phi$ is linear, and
as $U^*\AAA U=U_0^*\AAA U_0\oplus U_1^*\AAA U_1$,
$\phi$ is a bijection
(cf.\ Lang [6], Chapter XVI, Proposition 8.2).
Moreover, it is an algebra $*$-isomorphism, since
$U_i^*U_i=I$ for $i=0,1$ and hence
$U_0^{\otimes m}\otimes U_1^{\otimes n-m}$ is unitary.

Since a full block diagonalization of $\Sym^m(\oC^{2\times 2})$,
expressed in the standard basis of $\Sym^m(\oC^{2\times 2})$, is
known for any $m$ (Section \ref{9fe10b}.\ref{8fe10b}), and since
the tensor product of full block diagonalizations is again a full
block diagonalization, we readily obtain with $\phi$
a full block diagonalization of $\Sym^n(\AAA)$.
To use it in computations,
we need to describe it in terms of the standard basis of $\Sym^n(\AAA)$.
First we express $\phi$ in terms of the standard bases of
$\Sym^n(\AAA)$ and of $\Sym^m(\oC^{2\times 2})$ and
$\Sym^{n-m}(\oC^{2\times 2})$.

Let $\Lambda_T^n$ and $\Omega_{\lambda}$ be as in \rf{13fe10b} and
\rf{13fe10c}.
For $\lambda\in\Lambda_T^n$, define
\de{13fe10e}{
B_{\lambda}:=\sum_{\rho\in\Omega_{\lambda}}
\bigotimes_{i=1}^n B_{\rho(i)}.
}
Then $\{B_{\lambda}\mid\lambda\in\Lambda_T^n\}$ is a basis of $\Sym^n(\AAA)$.

We need the `Krawtchouk polynomial': for $n,k,t\in\oZ_+$,
\dez{
K_k^n(t):=\sum_{i=0}^k(-1)^i\kies{t}{i}\kies{n-t}{k-i}.
}
For later purposes we note here that for all $n,k,t$:
\dy{3fe10a}{
$K^n_{n-k}(t)=(-1)^tK^n_k(t)$.
}

For $\lambda\in\Lambda_T^n$, $\alpha\in\Lambda_P^m$,
$\beta\in\Lambda_P^{n-m}$, define
\dyy{23ja10a}{
\vartheta_{\lambda,\alpha,\beta}:=
\delta_{\lambda',\alpha+\beta}
\prod_{c\in P}
K^{\lambda'(c)}_{\lambda(c1)}(\beta(c)),
}
where for $\lambda\in\Lambda_T^n$, $\lambda'\in\Lambda_P^n$ is defined
by
\dez{
\lambda'(c):=\lambda(c0)+\lambda(c1)
}
for $c\in P$.

We now express $\phi$ in the standard bases \rf{13fe10e} and \rf{13fe10d}.

\propz{
For any $\lambda\in\Lambda_T^n$,
\dyy{15ja10b}{
\phi(B_{\lambda})
=
\bigoplus_{m=0}^n
\sum_{\alpha\in\Lambda_P^m,\beta\in\Lambda_P^{n-m}}
\vartheta_{\lambda,\alpha,\beta}
M_{\alpha}\otimes M_{\beta}.
}
}

\pf
By \rf{15ja10c}, the $m$-th component of $\phi(B_{\lambda})$
is equal to
\dyy{23ja10b}{
\sum_{\rho\in\Omega_{\lambda}}
\Big(
\bigotimes_{i=1}^m
E_{\rho_1(i),\rho_2(i)}
\Big)
\otimes
\Big(
\bigotimes_{i=m+1}^n
(-1)^{\rho_3(i)}
E_{\rho_1(i),\rho_2(i)}
\Big)
=
\sum_{\mu\in\Lambda_T^m,\nu\in\Lambda_T^{n-m}
\atop \mu+\nu=\lambda}
\Big(
\sum_{\sigma\in\Omega_{\mu}}
\bigotimes_{i=1}^m
E_{\sigma_1(i),\sigma_2(i)}
\Big)
\otimes
\Big(
\sum_{\tau\in\Omega_{\nu}}
\bigotimes_{i=1}^{n-m}
(-1)^{\tau_3(i)}
E_{\tau_1(i),\tau_2(i)}
\Big)
=
\sum_{\mu\in\Lambda_T^m,\nu\in\Lambda_T^{n-m}
\atop \mu+\nu=\lambda}
\Big(
\prod_{c\in P}\kies{\mu'(c)}{\mu(c1)}
\Big)
M_{\mu'}
\otimes
\Big(
\prod_{c\in P}(-1)^{\nu(c1)}\kies{\nu'(c)}{\nu(c1)}
\Big)
M_{\nu'}.
}
If we sum over $\alpha:=\mu'$ and $\beta:=\nu'$, we can next,
for each $c\in P$, sum over $j$ and set $\nu(c1):=j$, and
$\mu(c1):=\lambda(c1)-j$.
In this way we get that the last expression in \rf{23ja10b} is equal to
\dyyz{
\sum_{\alpha\in\Lambda_P^m,\beta\in\Lambda_P^{n-m}
\atop \alpha+\beta=\lambda'}
\Big(
\prod_{c\in P}
\sum_{j=0}^{\lambda(c1)}
(-1)^j
\kies{\alpha(c)}{\lambda(c1)-j}
\kies{\beta(c)}{j}
\Big)
M_{\alpha}\otimes M_{\beta}
=
\sum_{\alpha\in\Lambda_P^m,\beta\in\Lambda_P^{n-m}}
\vartheta_{\lambda,\alpha,\beta}
M_{\alpha}\otimes M_{\beta}.
\bxx
}

This describes the algebra isomorphism $\phi$ in \rf{22ja10b}.
With the results given in Section \ref{9fe10b}.\ref{8fe10b}
it implies a full block diagonalization
\dyy{3fe10b}{
\psi:
\Sym^n(\AAA)\to
\bigoplus_{m=0}^n
\bigoplus_{k=0}^{\lfloor\frac12m\rfloor}
\bigoplus_{l=0}^{\lfloor\frac12(n-m)\rfloor}
\oC^{[k,m-k]\times[k,m-k]}
\otimes
\oC^{[l,n-m-l]\times[l,n-m-l]},
}
described by
\dyy{22ja10c}{
\psi(B_{\lambda})=
\bigoplus_{m=0}^n
\bigoplus_{k=0}^{\lfloor\frac12m\rfloor}
\bigoplus_{l=0}^{\lfloor\frac12(n-m)\rfloor}
\psi_{m,k,l}(B_{\lambda})
}
where
\dyy{17me10c}{
\psi_{m,k,l}(B_{\lambda}):=
\sum_{\alpha\in\Lambda_P^m,\beta\in\Lambda_P^{n-m}}
\vartheta_{\lambda,\alpha,\beta}\Gamma_{\alpha,k}\otimes\Gamma_{\beta,l}
}
for $\lambda\in\Lambda_T^n$.
Inserting \rf{10fe10a} and \rf{22ja10a} in \rf{17me10c}
makes the block diagonalization explicit, and it can readily
be programmed.
Note that $\alpha,\beta$ in the summation can be restricted to those with
$\alpha+\beta=\lambda'$.
Note also that at most one entry of the matrix
$\Gamma_{\alpha,k}\otimes\Gamma_{\beta,l}$ is nonzero.

\subsectz{Deleting distances}

For $m,k,l$, we will use the natural isomorphism
\de{18me10a}{
\oC^{([k,m-k]\times[l,n-m-l])\times([k,m-k]\times[l,n-m-l])}
\cong
\oC^{[k,m-k]\times[k,m-k]}
\otimes
\oC^{[l,n-m-l]\times[l,n-m-l]}.
}

\prop{24fe10b}{
Let $D\subseteq[0,n]$.
Then
the linear hull of
\dez{
\{\psi_{m,k,l}(B_{\lambda})\mid\lambda\in\Lambda_T^n,
i_{\lambda'},j_{\lambda'}\in D\}
}
is equal to the subspace $\oC^{F\times F}$ of \rf{18me10a},
where
\dyyz{
F:=\{(i,i')\in[k,m-k]\times[l,n-m-l]\mid i+i'\in D\}.
}
}

\pf
For any $\lambda\in\Lambda_T^n$,
if $\psi_{m,k,l}(B_{\lambda})_{(i,i'),(j,j')}$ is nonzero, then
$i+i'=i_{\lambda'}$ and $j+j'=j_{\lambda'}$.
This follows from \rf{22ja10c} and
from the definition of the matrices $\Gamma_{\alpha,k}$
(cf.\ \rf{22ja10a}).

Hence, for any fixed $a,b\in\oZ_+$,
the linear hull of the $\psi_{m,k,l}(B_{\lambda})$ with
$i_{\lambda'}=a$ and $j_{\lambda'}=b$ is equal to the
the set of matrices in \rf{18me10a}
that are nonzero only in
positions $(i,i'),(j,j')$ with $i+i'=a$ and $j+j'=b$.
\bx 

So if distances are restricted to $D\subseteq[0,n]$, we can reduce
the block diagonalization to those rows and columns with index
in $F$.

\subsectz{Unordered pairs}

We now go over from ordered pairs to unordered pairs.
First, let $\NN'_2:=\NN_2\setminus\{\emptyset\}$, and
consider $\End_G(\oC^{\NN'_2})$.
Let $\tau$ be the permutation of $N^2$ swapping $(c,d)$ and $(d,c)$ in $N^2$.
Let $Q_{\tau}$ be the corresponding permutation matrix in $\oC^{N^2\times N^2}$.
Note that $N^2$ corresponds to the set of row indices of the matrices
$B_{\lambda}$ (cf.\ \rf{6ma10a}).

For any $\lambda\in\Lambda_T^n$,  let $\widetilde\lambda\in\Lambda_T^n$
be given by $\widetilde\lambda(t_1,t_2,t_3):=\lambda(t_1,t_2,t_3+t_1)$
for $t\in T$.
So for any $c\in P$, $\widetilde\lambda(c1)=\lambda(c1)$ if $c_1=0$
and $\widetilde\lambda(c1)=\lambda'(c)-\lambda(c1)$ if $c_1=1$.
Then $B_{\widetilde\lambda}=Q_{\tau}B_{\lambda}$.

Now \rf{3fe10a} gives that
for any $m$ and $\alpha\in\Lambda_P^m$, $\beta\in\Lambda_P^{n-m}$ one
has
\dez{
\vartheta_{\widetilde\lambda,\alpha,\beta}=
(-1)^{i_{\beta}}\vartheta_{\lambda,\alpha,\beta}.
}
This implies that the matrix
$\psi_{m,k,l}(B_{\lambda}+B_{\widetilde\lambda})$ has only
0's in rows whose index $(i,i')$ has $i'$ odd.
Similarly, the matrix $\psi_{m,k,l}(B_{\lambda}-B_{\widetilde\lambda})$ 
has only 0's in rows whose index $(i,i')$ has $i'$ even.
As the space of matrices invariant under permuting the rows
by $\tau$ is spanned by the matrices $B_{\lambda}+B_{\widetilde\lambda}$,
this space corresponds under $\psi_{m,k,l}$ to those matrices that
have 0's in rows whose index $(i,i')$ has $i'$ odd.

A similar argument holds for columns if we consider
$\hat\lambda(t_1,t_2,t_3):=\lambda(t_1,t_2,t_3+t_2)$.
Hence the image of $\psi_{m,k,l}$ of the set of
elements of $\Sym^n(\AAA)$ that are invariant under the operations
$\lambda\mapsto\widetilde\lambda$ and $\lambda\mapsto\hat\lambda$
is precisely equal to
\de{12me10a}{
\oC^{[k,m-k]\times [k,m-k]}\otimes\oC^{[l,n-m-l]_{\text{even}}\times [l,n-m-l]_{\text{even}}}.
}

\subsectz{Adding $\emptyset$}

So far we have a full block decomposition of $\End_G(\oC^{\NN'_2})$,
We need to incorporate $\emptyset$ in it.
It is a basic fact from representation theory that if
$V_1,\ldots,V_t$ is the canonical decomposition of $\oC^{\NN'_2}$
into isotypic components (cf.\ Serre [18]),
then $\End_G(\oC^{\NN'_2})=\bigoplus_{i=1}^t\End_G(V_i)$, and
each $\End_G(V_i)$ is $*$-isomorphic to a full matrix algebra.

We can assume that $V_1$ is the set of $H$-invariant elements
of $\oC^{\NN'_2}$.
Hence, as $\emptyset$ is $G$-invariant, $V_1':=\oC^{\emptyset}\oplus V_1$
is the set of $G$-invariant elements of $\oC^{\NN_2}$.
One may check that the block indexed by $(m,k,l)=(n,0,0)$ corresponds
to $V_1$.
So replacing block $(n,0,0)$ by $\End_G(V'_1)$ gives a full block
diagonalization of $\End_G(\oC^{\NN_2})$.
Note that $\End_G(V'_1)=\End(V'_1)$.

We can easily determine a basis for $V'_1$, namely the set
of characteristic vectors of the $G$-orbits of $\NN_2$.
Then for any $B\in \End_G(\oC^{\NN_2})$, we can directly calculate its
projection in $\End(V'_1)$.
This gives the required new component of the full block diagonalization.

\sect{12me10b}{Restriction to even words}

We can obtain a further reduction by restriction to the collection
$E$ of words in $\{0,1\}^n$ of even weight.
(The {\em weight} of a word is the number of 1's in it.)
By a parity check argument one knows that for even $d$ the bound
$A(n,d)$ is attained by a code $C\subseteq E$.
A similar phenomenon applies to $A_k(n,d)$:

\prop{24fe10a}{
For even $d\geq 2$,
the maximum value in \rf{2fe10a} does not change if $x(S)$ is required to
be zero if $S\not\subseteq E$.
}

\pf
Let $\varepsilon:N\to E$ be defined by $\varepsilon(w)=w$ if $w$ has even
weight and $\varepsilon(w)=w+e_n$ if $w$ has odd weight.
Here $e_n$ is the $n$-th unit basis vector, and addition is modulo 2.
If $d$ is even, then for all $v,w\in N$: $d_H(v,w)\geq d$ if and only if
$d_H(\varepsilon(v),\varepsilon(w))\geq d$.
Now $\varepsilon$ induces a projection $p:\oR^{\NN}\to\oR^{\EE}$,
where $\EE$ is the collection of codes in $\NN$ with all words having
even weight.

One easily checks that if $M_S(x)$ is positive semidefinite for all $S$,
then $M_S(p(x))$ is positive semidefinite for all $S$.
Moreover, $\sum_{v\in N}p(x)(\{v\})=\sum_{v\in N}x(\{v\})$.
\bx

This implies that restricting $x$ to be nonzero only on subsets $S$ of $E$
does not change the value of the upper bound.
However, it gives a computational reduction.
This can be obtained by using Proposition \ref{24fe10b} and by observing
that the restriction amounts to an invariance under an action of $S_2$,
for which we can use Section \ref{9fe10b}.\ref{8fe10a}.
The latter essentially implies that in \rf{12me10a} we can restrict the
left hand side factor to rows and columns with index in
$[k,m-k]_{\text{even}}$.
As it means a reduction of the program size by only a linear factor,
we leave the details to the reader.

\bigskip
\noindent
{\em Acknowledgement.}
We thank Niels Oosterling for very helpful comments on the method.

\section*{References}\label{REF}
{\small
\begin{itemize}{}{
\setlength{\labelwidth}{4mm}
\setlength{\parsep}{0mm}
\setlength{\itemsep}{1mm}
\setlength{\leftmargin}{5mm}
\setlength{\labelsep}{1mm}
}
\item[\mbox{\rm[1]}] E. Agrell, 
Bounds for unrestricted binary codes,\\
\url{http://webfiles.portal.chalmers.se/s2/research/kit/bounds/unr.html}

\item[\mbox{\rm[2]}] A.E. Brouwer, 
Table of general binary codes,\\
\url{http://www.win.tue.nl/~aeb/codes/binary-1.html}

\item[\mbox{\rm[3]}] P. Delsarte, 
{\em An Algebraic Approach to the Association Schemes of Coding Theory}
[Philips Research Reports Supplements 1973 No. 10],
Philips Research Laboratories, Eindhoven, 1973.

\item[\mbox{\rm[4]}] C.F. Dunkl, 
A Krawtchouk polynomial addition theorem and wreath product of
symmetric groups,
{\em Indiana University Mathematics Journal} 25 (1976) 335--358.

\item[\mbox{\rm[5]}] M. Gr\"otschel, L. Lov\'asz, A. Schrijver, 
{\em Geometric Algorithms and Combinatorial Optimization},
Springer, Berlin, 1988.

\item[\mbox{\rm[6]}] S. Lang, 
{\em Algebra --- Revised Third Edition},
Springer, New York, 2002.

\item[\mbox{\rm[7]}] J.B. Lasserre, 
An explicit equivalent positive semidefinite program for
nonlinear 0-1 programs,
{\em {SIAM} Journal on Optimization} 12 (2002) 756--769.

\item[\mbox{\rm[8]}] M. Laurent, 
Strengthened semidefinite bounds for codes,
{\em Mathematical Programming, Series B} 109 (2007) 239--261.

\item[\mbox{\rm[9]}] L. Lov\'asz, 
On the Shannon capacity of a graph,
{\em {IEEE} Transactions on Information Theory} {IT}-25 (1979) 1--7.

\item[\mbox{\rm[10]}] L. Lov\'asz, A. Schrijver, 
Cones of matrices and set-functions and 0--1 optimization,
{\em {SIAM} Journal on Optimization} 1 (1991) 166--190.

\item[\mbox{\rm[11]}] F.J. MacWilliams, N.J.A. Sloane, 
{\em The Theory of Error-Correcting Codes},
North-Holland, Amsterdam, 1977.

\item[\mbox{\rm[12]}] R.J. McEliece, E.R. Rodemich, H.C. Rumsey, Jr, 
The Lov\'asz bound and some generalizations,
{\em Journal of Combinatorics, Information \& System Sciences}
3 (1978) 134--152.

\item[\mbox{\rm[13]}] H.D. Mittelmann, 
An independent benchmarking of {SDP} and {SOCP} solvers,
{\em Mathematical Programming} 95 (2003) 407--430.

\item[\mbox{\rm[15]}] M. Nakata, B.J. Braams, K. Fujisawa, M. Fukuda, J.K. Percus, M. Yamashita, Z. Zhao, 
Variational calculation of second-order reduced density matrices by
strong N-representability conditions and an accurate semidefinite
programming solver,
{\em Journal of Chemical Physics} 128, 16 (2008) 164113.

\item[\mbox{\rm[14]}] NEOS Server for Optimization, 
\url{http://www-neos.mcs.anl.gov/}

\item[\mbox{\rm[16]}] A. Schrijver, 
 A comparison of the Delsarte and Lov\'asz bounds,
{\em {IEEE} Transactions on Information Theory} {IT}-25 (1979) 425--429.

\item[\mbox{\rm[17]}] A. Schrijver, 
New code upper bounds from the Terwilliger algebra and semidefinite
programming,
{\em {IEEE} Transactions on Information Theory} 51 (2005) 2859--2866.

\item[\mbox{\rm[18]}] J.-P. Serre, 
{\em Linear Representations of Finite Groups},
Springer, New York, 1977.

\item[\mbox{\rm[19]}] H.D. Sherali, W.P. Adams, 
A hierarchy of relaxations between the continuous and convex hull
representations for zero-one programming problems,
{\em {SIAM} Journal on Discrete Mathematics} 3 (1990) 411--430.

\item[\mbox{\rm[20]}] F. Vallentin, 
Symmetry in semidefinite programs,
{\em Linear Algebra and Its Applications} 430 (2009) 360--369.

\end{itemize}
}

\end{document}